\documentclass{amsart}

\usepackage{graphicx}

\newtheorem{theorem}{Theorem}[section]
\newtheorem{lemma}[theorem]{Lemma}

\newtheorem{corollary}[theorem]{Corollary}
\newtheorem*{Theorem}{Theorem}

\theoremstyle{definition}
\newtheorem{definition}[theorem]{Definition}

\theoremstyle{remark}
\newtheorem{remark}[theorem]{Remark}

\newcommand{\R}{\mathbb{R}}
\newcommand{\N}{\mathbb{N}}
\newcommand{\C}{\mathbb{C}}
\newcommand{\Z}{\mathbb{Z}}
\newcommand{\D}{\mathbb{D}}

\def\Jset{{Julia set }}
\newcommand{\sm}{\setminus}
\newcommand{\ovl}[1]{\overline{#1}}
\newcommand{\eps}{\varepsilon}
\newcommand{\Circle}{{\mathbb S}^1}
\newcommand{\I}{\mathcal{I}}

\usepackage{color}

\newcommand{\hide}[1]{}

\begin{document}

\title{Hausdorff Dimension and Biaccessibility for polynomial Julia Sets}

\author{Philipp Meerkamp}
\address{Department of Mathematics, Malott Hall, Cornell University,
Ithaca, NY 14853-4201}
\email{pmeerkamp@math.cornell.edu}

\author{Dierk Schleicher}
\address{School of Engineering and Science, Jacobs University,
Postfach 750~561, D-28725 Bremen, Germany}
\email{dierk@jacobs-university.de}



\begin{abstract}
We investigate the set of biaccessible points for connected
polynomial Julia sets of arbitrary degrees $d\geq 2$. We prove that
the Hausdorff dimension of the set of external angles corresponding
to biaccessible points is less than $1$, unless the Julia set is an interval.
This strengthens theorems of Stanislav Smirnov and Anna Zdunik: they 
proved that the same set of external angles has zero $1$-dimensional 
measure.
\end{abstract}

\maketitle

\section{Introduction}
\label{sec:Intro}

The filled-in Julia set $K$ of a polynomial $p\colon\C\to\C$, defined 
as the set of points in $\C$ with bounded orbits, is often an 
interesting set with rich topological, combinatorial or geometric 
properties. In many cases, such a set is a \emph{dendrite}: a compact 
connected locally connected set that does not disconnect the plane. 
In some sense, such a set can often be viewed as an infinite tree. 
One way to ask our main question is, which proportion of this tree 
consists of ``endpoints'', and which proportion consists of 
``non-endpoints''? For finite (non-degenerate) trees, there are 
finitely many endpoints and a continuum of non-endpoints (on the 
arc); but what if the Julia set is a dendrite? And what in more 
general cases of filled-in Julia sets? And how are endpoints defined 
in general?

One way to make this definition precise is to say that a point $z$ in 
a tree is an endpoint unless there are points $x,y$ in the tree, 
different from $x$, such that $z$ is on an injective path connecting 
$x$ to $y$. This definition works best for path-connected Julia sets. 
The definition that we use is the following: \emph{a point $x\in K$ 
is an endpoint unless there are two curves in $\C\sm K$, not 
homotopic to each other, that connect $\infty$ to $x$}. Such points 
$x$ are also called \emph{biaccessible}; an equivalent definition is 
that two different dynamic rays land at $x$ (see the next section). 
Our main result is that in almost all cases, most points in $K$ (in a 
very strong sense) are endpoints, unless $K$ is an interval (a 
straight line segment).

\begin{Theorem}[Hausdorff dimension of biaccessible angles]
\label{Thm:Main}
Let $p$ be a polynomial of degree at least $2$ with connected Julia 
set $J=\partial K$. Then the
biaccessible points have external angles in a set of Hausdorff
dimension less than 1, except when $J$ is an interval.
\end{Theorem}

This result does not have any topological hypotheses on the Julia 
set, other than that it be connected. It need not be locally 
connected or path connected, or uniquely path connected. One might 
thus wonder whether a similar statement might be true for planar 
dendrites (or even continua) with certain conformal self-similarity 
properties, whether arising in complex dynamics or elsewhere.

Our result extends known theorems by several people. The fact that 
biaccessible points have external angles in a set of $1$-dimensional 
measure zero was shown by Smirnov \cite{Sm} and Zdunik \cite{Zd} (in 
other words, these points have zero harmonic measure), except if the 
Julia set is a straight line segment. Earlier, Zakeri \cite{Za} had 
shown this for quadratic polynomials with locally connected Julia 
sets.

The Hausdorff dimension of the set of external angles of biaccessible 
points has been investigated as well. Zakeri \cite{Za1} estimated 
this dimension for certain real quadratic polynomials, and Bruin and 
Schleicher \cite[Section~14]{BS}, \cite{BSnew} gave estimates for all 
complex quadratic polynomials, as well as for certain subsets of 
parameter space (the Mandelbrot set). Radu \cite{RR}, in his Bachelor 
thesis, proved the same result as ours for the case of connected 
polynomial Julia sets satisfying certain additional hypotheses, 
including local connectivity of the Julia set and further assumptions 
on the critical values.

\remark
There are also some known lower bounds on the Hausdorff dimension of 
external angles of biaccessible points. The dimension is clearly zero 
for  Julia sets with no or only countably many biaccessible points. 
Bruin and Schleicher \cite[Section~14]{BS}, \cite{BSnew} also showed 
that the dimension is zero for quadratic ``Feigenbaum'' polynomials: 
these are limits of Julia sets with only countably many biaccessible 
points. For all other quadratic polynomials with connected Julia 
sets, the dimension is strictly positive. It is natural to ask how 
these results extend to polynomials of higher degrees, and whether 
the Hausdorff measure of the set of biaccessible
angles in the ``right'' dimension is finite and positive (at least 
when the dimension is strictly positive).

Bruin and Schleicher \cite[Section~14]{BS}, \cite{BSnew} also showed 
that there can be no uniform upper bound on the Hausdorff dimension 
of biaccessible angles, even for fixed degrees: there are quadratic 
polynomials close to $z\mapsto z^2-2$ where the  biaccessible points 
have Hausdorff dimension arbitrarily close to $1$.

Further properties of biaccessible points in polynomial Julia sets 
have been studied by Zakeri and Schleicher \cite{SZ, Za}.

\remark
The concept of biaccessibility is a topological one, defined in terms 
of homotopy classes of curves outside of the Julia set. However, all 
our arguments are combinatorial and would allow to restate the result 
in a combinatorial way (that would actually be slightly stronger 
because not all dynamic rays must land). We will briefly discuss this 
in Section~\ref{Sec:Laminations}.

\medskip\noindent
\emph{Acknowledgement}. We would like to thank Remus Radu, Saeed 
Zakeri and especially Xavier Buff for useful discussions on this 
topic and on an earlier draft of this paper.


\section{Background}

Let $p: {\C} \rightarrow
{\C}$ be a polynomial of degree $d\geq 2$, which we may as
well assume to be monic.
The {\it filled \Jset}of $p$ is the set of all points in $\C$ with bounded
forward orbit under $p$,
and the {\it \Jset}of $p$ is $J:=\partial K$.  It satisfies
$J=p(J)=p^{-1}(J)$ \cite[Lemma~4.3]{M}.

The sets $K$ and $J$ are nonempty compact subsets of $\C$.
Moreover, $K$ is full, i.e., $\C\setminus\ K$ is connected. In this paper, we
assume that $J$ (or equivalently $K$) is connected.
In this case, there is a Riemann map
\begin{equation}
\psi: \overline{\C} \setminus \overline{\D} \rightarrow
\overline{\C} \setminus K \;.
\end{equation}
When $p$ is monic, then $\psi$ is unique by requiring
$\psi(\infty)=\infty$ and $\lim_{z\to\infty} \psi(z)/z\to 1$. This
Riemann map satisfies
\begin{equation}
\psi(z^d) = p(\psi(z)) \;;
\label{Eq:Conjugation}
\end{equation}
i.e., it conjugates $p$ on its basin of $\infty$ to $z\mapsto z^d$.

For $t\in\Circle:=\R/\Z$, the image of the radial line $\{\psi(re^{2 
\pi it}) : r>1 \}$ is
called the {\it dynamic ray at (external) angle $t$} and denoted $R_t
\subset \overline{\C} \setminus K$; it satisfies $p(R_t) = R_{dt}$.

Consider the radial limit
\begin{equation}\label{Limit}
\gamma(t)\ =\ \lim_{r\searrow 1}\psi(re^{2\pi i t}).
\end{equation}
It need not exist for all $t\in\Circle$, but it is well known to exist
for almost all $t$ \cite[Theorem~18.2]{M}. If this limit exists, one
says that the ray $R_t$ {\it lands} at the
point $\gamma(t)\in J$.

For all angles $t\in \Circle$ for which $R_t$ lands, the ray $R_{dt}$ 
lands as well, and
\begin{equation}\label{conjugacy}
\gamma(dt)\ =\ p(\gamma(t)) \;.
\end{equation}
Thus $\gamma$ is a semiconjugation of multiplication by $d$ on 
$\Circle$ to the action of $p$ on the
\Jset (restricted to those angles the rays of which land, and the 
corresponding landing points).  In the particular
case when the Julia set $J$ is locally connected, the map $\gamma :
\Circle \rightarrow J$ is defined everywhere, and it is a continuous
surjection $\gamma\colon \Circle\to J$  \cite[Theorem~18.3]{M}.

We will denote distance on $\Circle$ by $\tau$ (normalized so as to
inherit the metric locally from $\R$, and always measuring along the
shorter of the two arcs connecting two points in $\Circle$). For an
interval $I\subset \Circle$, let $\tau(I)$ denote its length.

A point $z\in J$ is called {\it accessible} if $z$ is the landing
point of a dynamic ray (by Lindel\"of's theorem, this is equivalent
to the existence of a curve in $\C\setminus K$ converging to $z$).
The point $z$ is called {\it biaccessible} if it is the landing point
of at least two rays. If two dynamic rays $R_t$ and $R_{t'}$ land together,
we call the external angles $t$ and $t'$ {\it biaccessible angles},
and we call $(t,t')$ a {\it biaccessible angle pair}. Let 
$\Lambda\subset\Circle$ be the  set of all angles in biaccessible 
angle pairs. A \emph{ray pair} is a set of two dynamic rays that land 
at a common point (so that their external angles form a biaccessible 
angle pair).

As usual, a  {\it critical point} of $p$ is a point $z$ with
$p'(z)=0$, and $p(z)$ is the corresponding {\it critical value}.

\begin{remark}\label{Remark}
Note that $z$ is biaccessible if and only if $p(z)$ is biaccessible, unless $z$
is a critical point for $p$. This also means that $t$ is a
biaccessible angle if and only if so is $dt$, except when $R_t$ lands at a
critical point.
\end{remark}

\begin{remark}
\label{Rem:PolyAccessible}
It may also happen that three or more rays land at a single point. We 
call such a point ``poly-accessible''. It is well known and easy to 
see that such points always form a countable set 
\cite[Proposition~2.18]{Po}
(for every $\eps>0$, there can only be finitely many landing points 
of three rays the angles of which have mutual distance at least 
$\eps$). If three or more rays land at the same point that is not 
eventually periodic and not eventually critical, then this landing 
point is called a ``wandering triangle'' (or, more generally, a 
``wandering polygon''), and the number of rays at such points, as 
well as the number of such orbits, satisfy certain bounds depending 
on the degree of $d$, in particular these numbers are finite (compare 
\cite{Kiwi,BL,BCLOS}). If three or more rays land at a periodic 
point, then either the landing point is repelling or parabolic and 
the number of rays is finite, or the landing point is a Cremer point 
(in which case it is not known whether any rays may land at all; 
compare \cite{SZ}).
In total, there can thus be only countably many poly-accessible 
points, and exceptfor Cremer points the total number of rays involved 
is countable.
\end{remark}


\section{Endpoints of the Julia Set}
\label{sec:Arcs}

\begin{definition}[J-Endpoint]
\label{def:ends}
A point $t\in \Circle$ will be called a {\it J-endpoint} of a Julia 
set if there exists a
sequence $(t_n,t'_n)$ of biaccessible angle pairs such that
\begin{displaymath}
t_n \rightarrow t, \hspace{20pt}
t'_n \rightarrow t,
\end{displaymath}
and such that for all large $n$, the point $t$ lies in the shorter
arc of $\Circle$ connecting $t_n$ and $t'_n$.
\end{definition}

\begin{lemma}[Trichotomy]
\label{Lem:Trichotomy}
For every connected polynomial Julia set $J$, exactly one of the
following three cases holds:
\begin{itemize}
\item
$J$ has no biaccessible points and $\Lambda$ is empty;
\item
$J$ is an interval and $\Lambda$ is $\Circle$ minus two points;
\item
$J$ has at least three J-endpoints and $\Lambda$ is dense.
\end{itemize}
\end{lemma}
\begin{proof}
These three cases are clearly mutually exclusive. If $\Lambda$ is
non-empty, then the preimage of any biaccessible point in $J$ is a
biaccessible point, and it follows easily that $\Lambda$ is infinite
and dense.

Now we show that if $\Lambda$ is non-empty, then $J$ has at least two
J-endpoints. Consider a biaccessible angle pair $(t,t')$. The
angles $t$ and $t'$ separate the circle $\Circle$ into two open intervals, say
$I$ and $I'$. Since $\Lambda$ is
dense, there is a biaccessible angle pair $(t_1,t_1')$ with
$t_1,t'_1\in I$. Let $I_1\subset I$ be the
interval bounded by $t_1$ and $t_1'$. It has strictly
smaller length than $I$; in fact, choosing $t_1$ in the middle
third of $I$, we can assure that the length of $I_1$ is at most
$2/3$ of the length of $I$. Iterating this argument, we find a
biaccessible
angle pair in arbitrarily small intervals inside $I$ and hence at
least one J-endpoint in $I$ (i.e., this J-endpoint is in the 
impression of a ray with angle in $I$). In a similar way we find at
least one J-endpoint in $I'$. Thus the Julia set $J$ has at least two
J-endpoints.

If $J$ has exactly two J-endpoints, then $J$ is an interval by
\cite[Lemma~2]{Zd} (and $p$ is a Chebyshev polynomial, up to sign);
see also \cite{Zd1}. (The idea is this: if a critical value is not a
J-endpoint, then the corresponding critical point is a branch point of
$J$ and creates extra J-endpoints; moreover, J-endpoints always
map to J-endpoints. This implies that if $J$ has only two J-endpoints,
then it is postcritically finite and has only two critical values. By 
conjugation, we may suppose that these two critical values are real.
The polynomial $p$ has a Hubbard tree without branch points, and all 
critical values are endpoints of this tree. The Hubbard tree is thus 
backwards invariant, so it equals the Julia set.)
\end{proof}

For us, the only interesting case is when the Julia set has at least
three J-endpoints. We will assume this from now on.

\section{Narrow Preimages of Ray Pairs}
\label{Sec:Narrow}

Let $(t,t')$ be a biaccessible angle pair and $(t_1,t'_1)$ one of its 
preimage biaccessible angle pairs, i.e., a biaccessible angle pair 
with $dt_1=t$ and $dt'_1=t'$. We call the preimage $(t_1,t'_1)$ 
\emph{narrow} if $\tau(t_1,t'_1)=\tau(t,t')/d$. An (iterated) 
preimage biaccessible angle pair $(t_n,t'_n)$ is called \emph{narrow 
of generation $n$} if $(t_n,t'_n)$ is obtained from $(t,t')$ by 
taking $n$ generations of preimages, and all intermediate preimages 
are narrow (so that $\tau(t_n,t'_n)=\tau(t,t')/d^n)$\,).

Our assumption that the Julia set has at least three J-endpoints 
implies that there are three biaccessible angle pairs, say $(a,a')$, 
$(b,b')$ and $(c,c')$, so that none of them separates the other two: 
we may assume that the cyclic order of these six angles is 
$a,a',b,b',c,c'$. To simplify the reasoning, suppose that the longer 
of the two intervals of $\Circle\sm\{a,a'\}$ contains 
$\{b,b',c,c'\}$, and similarly for $\Circle\sm\{b,b'\}$ and 
$\Circle\sm\{c,c'\}$. We will call biaccessible angle pairs with this property 
\emph{un-nested}. Each Julia set with at least three J-endpoints 
clearly has three such angle pairs.

\begin{lemma}[Number of narrow preimage biaccessible angle pairs]
\label{Lem:NumberStronglyNarrow}
Consider a polynomial $p$ of degree $d\ge 2$ and suppose it has three 
un-nested biaccessible angle pairs. For generations $n\ge 0$, let 
$s_n$ be the combined number of narrow preimage biaccessible angle 
pairs of all three biaccessible angle pairs. Then $s_{n+1}\ge 
ds_n-2(d-1)$.
\end{lemma}
\begin{proof}
Denote the three given biaccessible angle pairs by $(a,a')$, $(b,b')$ 
and $(c,c')$.
The set $\C\sm \ovl{R_a\cup R_{a'}}$ consists of two components, say 
$B_a$ and $B'_a$. The assumption that the three given ray pairs are 
non-nested means that one of these components, say $B'_a$, contains 
the other two given ray pairs; define $B_b$ and $B_c$ analogously.

Suppose for simplicity that non of the ray pairs considered lands at 
a critical value.
The immediate preimage of the biaccessible angle pair $(a,a')$ 
consists of $d$ disjoint biaccessible angle pairs that may or many 
not be nested. The corresponding ray pairs disconnect $\C$ into $d+1$ 
open complementary domains $U_1,\dots, U_{d+1}$. Each of the $U_j$ 
maps under $p$ as a proper holomorphic map onto a component of $\C\sm 
\ovl{R_a\cup R_{a'}}$, so each $U_j$ has an associated mapping degree 
that equals the number of critical points on $U_j$ plus $1$ (always 
counting multiplicities). A domain $U_j$ is narrow (i.e., $U_j$ is 
bounded by a single ray pair, and this ray pair is narrow) if and 
only if $U_j$ maps conformally onto $B_a$; this implies that $U_j$ 
does not contain a critical point.

The case with the maximal count of narrow components $U_j$ is when a 
single domain $U_{d+1}$ contains all $d-1$ critical points and maps 
onto $B'_a$ with degree $d$, and all other domains $U_1$,\dots, $U_d$ 
map conformally onto $B_a$ and are narrow.

If a domain $U_j$ contains a single critical point, then its mapping 
degree is $2$; if $U_j$ maps to $B_a$, then the count of possible 
narrow domains $U_j$ decreases by $2$. Each further critical point in 
$U_j$ decreases this count by $1$. Therefore, the number of narrow 
components $U_j$ is at least $d$ minus twice the number of critical 
points the images of which are in $B_a$.

The same arguments apply to the angle pairs $(b,b')$ and $(c,c')$. 
Note that the domains $B_a$, $B_b$, and $B_c$ are disjoint, so each 
critical point that reduces the number of narrow components can count 
only for one of the three ray pairs. Therefore, the total number of 
narrow preimages of $(a,a')$, $(b,b')$ and $(c,c')$ is $s_1\ge 
3d-2(d-1)=d+2$. All $s_1$ ray pairs are narrow and non-nested. Note 
that initially we had $s_0=3$ angle pairs, and $s_1\ge ds_0-2(d-1)$ 
as claimed.

For the inductive step, the three initial angle pairs are replaced by 
$s_n$ non-nested ray pairs; for these, the argument can be repeated: 
each of these $s_n$ has at least $ds_n$ narrow preimages minus twice 
the number of critical points that map into the narrow complementary 
component of any ray pair, and again each critical point can reduce 
the count of narrow preimages for only one ray pair, and only by two. 
This yields the formula $s_{n+1}\ge ds_n-2(d-1)$ as claimed.

Finally, if some ray pair lands at a critical value, then some of its 
preimage ray pairs merge, and the statement remains true for an 
appropriate choice of rays to form pairs.
\end{proof}

Let $E$ be the set of external angles with the property that for each 
$t\in E$, there are infinitely many narrow biaccessible angle pairs 
$(t_k,t'_k)$ on the backwards orbit of $(a,a')$, $(b,b')$, or 
$(c,c')$ so that the interval $(t_k,t'_k)\subset\Circle$ contains 
$t$, and so that the lengths of these intervals tend to $0$. ($E$ 
stands for ``endpoints'' of the Julia set.)

\begin{lemma}
Angles in $E$ are not part of any biaccessible angle pair.
\end{lemma}
\begin{proof}
Suppose that $(t,t')$ is a biaccessible angle pair, but $t\in E$. 
Then there is a narrow biaccessible angle pair $(t_k,t'_k)$ with 
$\tau(t_k,t'_k)<\tau(t,t')$ andso that the interval $(t_k,t'_k)$ 
contains $t$; this interval thus cannot contain $t'$. This is a 
contradiction unless all the rays at angles $t_k,t'_k,t,t'$ land at 
the same point. This would imply that for all biaccessible angle 
pairs $(t_m,t'_m)$ of higher generation than $(t_k,t'_k)$ such that 
$t$ is contained in the interval $(t_m,t'_m)$, the rays $R_{t_m}$ and 
$R_{t'_m}$ must also land at the same point, so infinitely many rays 
would land at this point. We will now show that this is impossible.

Indeed, the set $\C\sm\ovl{R_{t_k}\cup R_{t'_k}}$ consists of two 
components; let $U_k$ be the component containing $R_t$, and define 
$U_m$ similarly. We may assume that both angle pairs are on the 
backwards orbit of the same angle pair $(a,a')$, $(b,b')$, or 
$(c,c')$, and indeed that $(t_m,t'_m)$ is on the backwards orbit of 
$(t_k,t'_k)$, say after $s>0$ iterations. Since all biaccessible 
angle pairs are narrow, this implies that $p^{\circ s}\colon U_m\to 
U_k\supset U_m$ is a conformal isomorphism with conformal inverse 
$q\colon U_k\to U_m$. Iterating $q$, it follows that $t$ is a 
periodic angle, so the landing point is repelling or parabolic, and 
only finitely many rays land there (compare 
Remark~\ref{Rem:PolyAccessible}).
\end{proof}

So all we need to do is prove that $\Circle\sm E$ has Hausdorff 
dimension less than $1$ (most angles correspond to J-endpoints). We 
will prove the following result.

\begin{lemma}
\label{Lem:DimensionLemma}
The Hausdorff dimension of $\Circle\sm E$ satisfies 
$\dim_H(\Circle\sm E)<\eta<1$, where $\eta$ depends only on $d$ and 
the lengths of $(a,a')$, $(b,b')$, and $(c,c')$. \end{lemma}

\begin{remark}
Let $\alpha$ be the minimum of the three lengths of $(a,a')$, 
$(b,b')$, $(c,c')$.
The quantity $\alpha$ is, in some sense, a measure of the ``size'' of 
branch points of the Julia set: in a locally connected Julia set, a 
branch point $q$ is the landing point of at least three dynamic rays, 
and its size can be defined as
\[
s(q)=\sup\{\delta>0\colon \mbox{three rays at angles $t_1$, $t_2$, 
$t_3$ land at $q$ with $\tau(t_i,t_j)>\delta$ for $i\neq j$}\} \;.
\]
\end{remark}
This definition coincides with the maximal possible value of $\alpha$ 
in the locally connected case, so $\alpha$ can indeed be seen as a 
measure of size of a branch point (a branch point is ``small'' if all 
the rays landing at it can be grouped into two sets so that the 
angles are contained in short intervals).

The proof of Lemma~\ref{Lem:DimensionLemma} will be given, in 
somewhat abstract form, in the next section.


\section{Hausdorff dimension}

Let $g\colon\Circle\to\Circle$ be multiplication of angles by $d\ge 2$.

\begin{lemma}
Let $I_1,I_2,I_3\subset\Circle$ be three disjoint open intervals, 
each of length $d^{-N}$ for some $N\in\N$, and each bounded by 
endpoints in $\Z/d^N$ (projected into $\R/Z$).
Let $\I_0:=\{I_1,I_2,I_3\}$. For $n\ge 1$, let $\I_n$ be the set of 
all $g$-preimages of all intervals in $\I_{n-1}$, except that 
$2(d-1)$ intervals are missing in $\I_n$. Then $|\I_n|=d^{n}+2$, and 
for every $m\in\N$, the set
\[
\Circle\sm\bigcup_{n\ge m} \bigcup_{I\in\I_n}I
\]
has Minkowski dimension at most $1-3/{d^{N}}{N\log d}<1$.
\end{lemma}
\begin{proof}
We have $|\I_0|=3=d^0+2$ and $|\I_{n+1}|=d|\I_n|-2(d-1)=d^{n+1}+2$ by 
induction.

The condition on the endpoints assures that whenever $n'>n$, each 
interval in $\I_{n'}$ is either disjoint from or contained in any 
interval in $\I_n$.

For $n\in\N$, set $A_n:=\bigcup_{i\in\I_n} I$, and let $c_n$ be the 
number of intervals of length $d^{-n}$ required to cover $\Circle\sm 
A_n$.

The set $\Circle\sm A_0$ consists of $3$ intervals of length less 
than $1$. The set $\Circle\sm g^{-n}(A_0)$ thus consists of $3d^n$ 
intervals of length less than $d^{-n}$, for each $n$ (in this 
argument, we do not delete the $2(d-1)$ intervals in each generation).

In general, we have $c_{n+1}\le dc_n+2(d-2)$: taking preimages under 
$g$ increases the required number of intervals by a factor of $d$ 
(the length of the intervals decreases by a factor $d$), and each 
removed interval in $\I_{n+1}$ might require an extra covering 
interval.

Choose $\beta>0$ so that $\gamma:=(1+\beta)^N(1-3d^{-N})<1$. There is 
an $M\in\N$ so that $c_{n+1}\le (1+\beta) d c_n$ for all $n\ge M$; 
possibly by enlarging $M$, suppose that $M$ is divisible by $N$.

The set $\Circle\sm A_{M}$ can be covered by some number $C:=c_{M}$ 
intervals of length $d^{-M}$. Then $\Circle\sm A_{M+N}$ can be 
covered by $c_{M+N}\le d^N (1+\beta)^NC$ intervals of length 
$d^{-(M+N)}$, and $\Circle\sm(A_{M+N}\cup A_0)$ can be covered by 
$(d^N-3) (1+\beta)^NC=\gamma d^N C$ intervals of the same length, 
because one in $d^N$ intervals is contained in each of the intervals 
in $A_0$.

This argument can be repeated: $\Circle\sm(A_{M+2N}\cup A_N)$ can be 
covered by $(1+\beta)^N\gamma d^{2N}C$ intervals of length 
$d^{-(M+2N)}$, and for $\Circle\sm(A_{M+2N}\cup A_N\cup A_0)$, at 
most $(1-3d^{-N})(1+\beta)^N\gamma d^{2N}C=\gamma^2 d^{2N}C$ 
intervals are required.

Let us focus on the case $m=0$ in the claim:
we are interested in the set $\Circle\sm\bigcup_{n\ge 0} A_n$,  and 
this is contained in
\[ \Circle\sm\bigcup_{k\ge 0} A_{kN} \subset 
\Circle\sm\left(A_{M+KN}\cup A_{(K-1)N}\cup A_{(K-2)N}\cup \dots\cup 
A_N\cup A_0 \right) \;,
\]
and this latter set can be covered by $\gamma^Kd^{KN}C$ intervals of 
length $d^{-(M+KN)}$.

For the Minkowski dimension, we get the upper bound
\[
\lim_{K\to\infty} \frac{\log (\gamma^Kd^{KN}C)}{\log d^{M+KN}} = 
\frac{\log(\gamma d^N)}{\log d^{N}} = 
1+\frac{\log(1-3d^{-N})+N\log(1+\beta)}{\log d^N} \;.
\]
Now $\beta$ can be chosen arbitrarily close to $0$. Since 
$\log(1-3d^{-N})<-3d^{-N}$, the dimension is bounded above by 
$1-3d^{-N}/N\log d$.

Now we treat the case $m>0$. The set $\Circle\sm \bigcup_{n\ge 0} 
A_{n}$ can be covered by $\gamma^Kd^{KN}C$ intervals of length 
$d^{-M+KN}$, for every $K\ge 0$, so $\Circle\sm\bigcup_{n\ge m}A_n$ 
can be covered by $(1+\beta)^m\gamma^Kd^{KN+m}C$ intervals of length 
$d^{-(M+KN+m)}$, and this leads to the same dimension.
\end{proof}

\begin{corollary}
Using the notation of the previous lemma, the set
\[
\Circle\sm \bigcap_{m\ge 0} \bigcup_{n\ge m}\bigcup_{I\in\I_n}
\]
has Hausdorff dimension at most
$\displaystyle 1-{3}/{d^NN\log d}$.
\end{corollary}
\begin{proof}
By the lemma, each set $B_m:=\Circle\sm \bigcup_{n\ge 
m}A_n=\Circle\sm\bigcup_{n\ge m}\bigcup_{I\in\I_n} I$ satisfies the 
same upper bound for Minkowski dimension, hence Hausdorff dimension. 
The set we are interested in is $\Circle\sm\bigcap_m\bigcup_{n\ge m} 
A_m =\bigcup_m(\Circle\sm\bigcup_{n\ge m} A_n) = \bigcup_m B_m$, and 
Hausdorff dimension is stable under countable unions.
\end{proof}

\begin{proofof}{Lemma~\ref{Lem:DimensionLemma}}
If the Julia set has three J-endpoints, then there are three 
biaccessible angle pairs $(a,a'), (b,b'), (c,c')$ so that the three 
intervals $(a,a')$, $(b,b')$, $(c,c')$ (as subsets of $\Circle$) are 
disjoint. Let $\alpha>0$ be the minimum of their lengths.
Their combined number of narrow preimages of generation $1$ is at 
least $3d-2(d-2)=d+2$ by Lemma~\ref{Lem:NumberStronglyNarrow}. If 
$d\ge 5$, then at least one of these three intervals, say $(a,a')$, 
has three narrow preimages of length $\alpha/d$.

Call these intervals $I_1$, $I_2$, $I_3$. The number of narrow 
preimages of further generations grows as in 
Lemma~\ref{Lem:NumberStronglyNarrow}. If we construct sets of 
intervals $\I_n$ as above  and set $A_n:=\bigcup_{I\in\I_n} I$, then 
we have $E=\bigcap_{m\ge 0}\bigcup_{n\ge m} A_n$.

Our intervals do not yet satisfy the condition on the form of the 
endpoints, so we cannot directly apply the corollary. Restricting 
$(a,a')$ to a subinterval of length at least $1/2d$ times the 
original length, we obtain an interval $I_0\subset(a,a')$ that is 
bounded by two numbers $k/d^N$ and $(k+1)/d^N$ for $k,N\in\N$. This 
yields smaller sets $A_n$ and a smaller set $E$, hence a larger set 
$\Circle\sm E$. The corollary applied to this larger set shows that 
$\Circle\sm E$ has Hausdorff dimension less than $1$.

If $d\in\{3,4\}$, then we might have to resort to narrow intervals of 
generation two, ofwhich there are at least $d^2+2\ge 11$ for all 
three intervals combined, and the argument proceeds as above (except 
that the dimension formula uses intervals of length $\alpha/d^2$, 
rather than $\alpha/d$).

Finally, if $d=2$, then we have to go one generation further.
\qed
\end{proofof}
\medskip

This also proves the theorem.

\section{A Remark on Laminations and Entropy}
\label{Sec:Laminations}

As mentioned earlier, our results can also be stated in terms of 
laminations of Julia set, as developed by Thurston 
\cite[Chapter~II]{Th}. In the complex unit disk $\D$, we identify the 
boundary $\partial\D=\Circle/\Z$ with the set of external angles. For 
a given connected polynomial Julia set $J$, we connect two angles 
$\alpha,\beta\in\partial\D$ by a geodesic of $\D$ whenever the two 
rays at angles $\alpha$ and $\beta$ land at a common (biaccessible) 
point. One usually uses hyperbolic geodesics for this purposes 
because this yields clearer pictures, even though the results are 
equivalent. Every biaccessible point is thus represented by an arc in 
$\D$; polyaccessible points are represented by polygons. The 
\emph{lamination} $J$ is the closure of all these arcs (if the Julia 
set is locally connected, then the set of arcs is closed to begin 
with), and each arc is called a \emph{leaf} in this lamination. A 
degenerate leaf is one that connects an angle to itself (it is a 
point).

\begin{figure}[tp]
\includegraphics[width=0.5\textwidth]{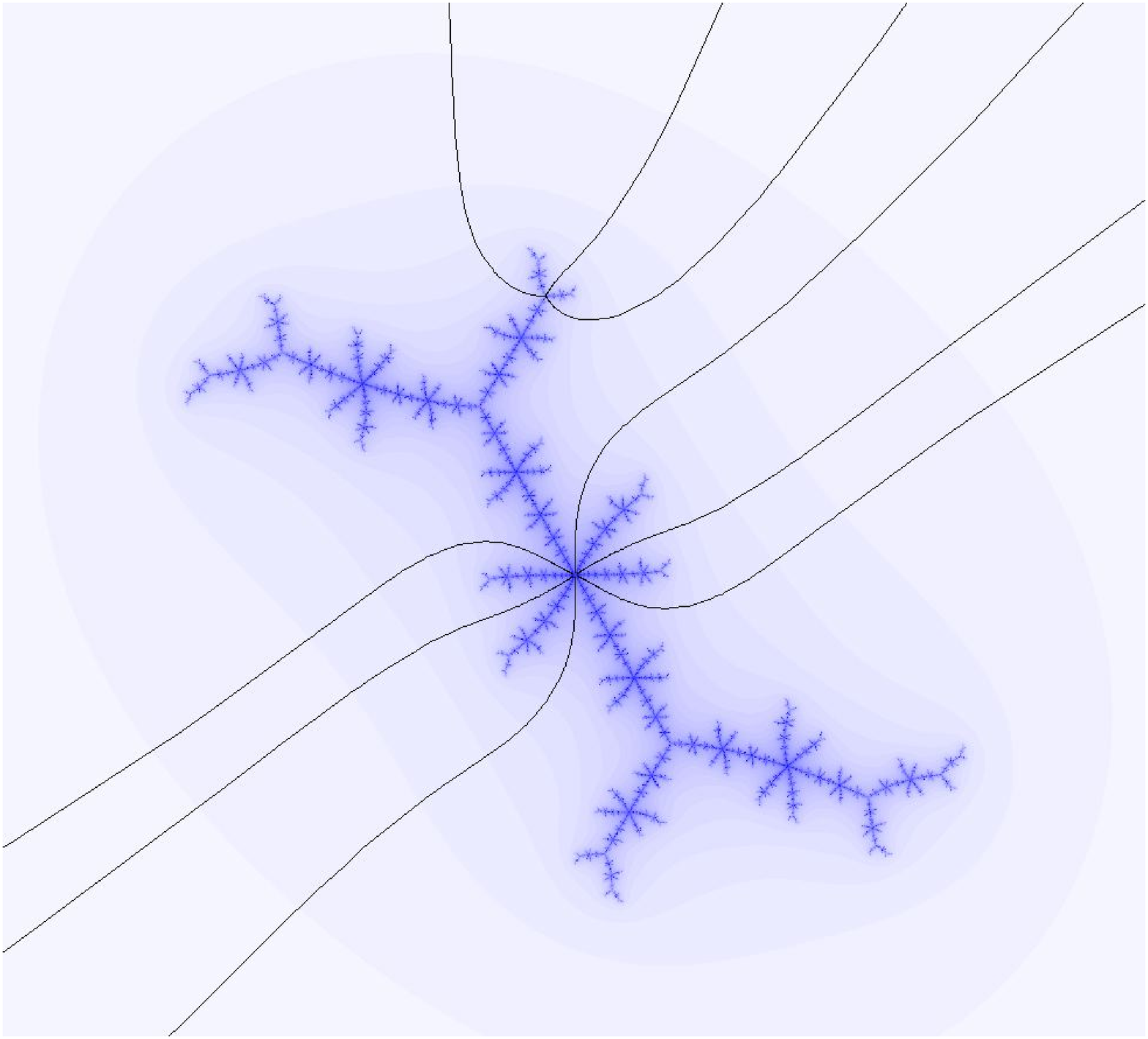}
\includegraphics[width=0.49\textwidth,trim=10 0 0 0,clip]{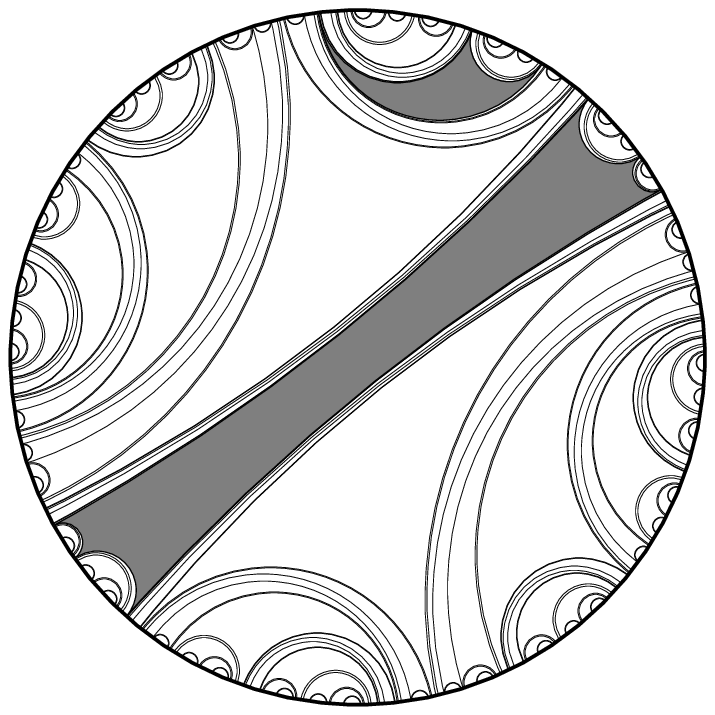}
\caption{Left: a quadratic polynomial with three rays landing at the 
critical point, and thus six rays landing at the critical value. 
Right: a lamination modeling this Julia set. The critical point 
corresponds to a hexagon, the critical value to a triangle (both 
shaded).
}
\label{FigGapPreperiodicPolygon}
\end{figure}

All our arguments, especially those in the key 
Section~\ref{Sec:Narrow}, are combinatorial in nature and can thus be 
stated entirely in terms of Julia set laminations, or more 
combinatorially in terms of abstract laminations that satisfy certain 
natural invariance properties as discussed in \cite[Chapter~II]{Th} 
(for the relation of laminations to Julia sets, see the appendix in 
\cite{Th}). We thus obtain an estimate about the Hausdorff dimension 
of the set endpoints of non-degenerate leaves in an invariant 
lamination. This setting relates to recent work of Thurston (as 
communicated in various talks and personal communication) as follows: 
if this dimension of endpoints is $\eta$, then the union of the 
leaves has dimension $1+\eta$. This dimension is related to the 
\emph{core entropy} of a postcritically finite Julia set: such a 
Julia set has a Hubbard tree (defined as a minimal invariant tree 
connecting all postcritical points), and the core entropy is the 
entropy of the polynomial restricted to the Hubbard tree. If $H$ 
denotes the core entropy, then $\eta=H/\log d$. Thurston's formula 
thus relates the topological concept of biaccessible points, and more 
precisely the geometric concept of the Hausdorff dimension of their 
angles, to the dynamic concept of entropy on the Hubbard tree.

\bibliographystyle{amsplain}

\begin{thebibliography}{99}

\bibitem{BCLOS} A. Blokh, C. Curry, G. Levin, L. Oversteegen, and D. 
Schleicher,
\emph{An extended Fatou-Shishikura inequality and wandering branch 
continua for polynomials}. Manuscript, submitted.

\bibitem{BL} A. Blokh, G. Levin,
\emph{An inequality for laminations, Julia sets and ``growing trees''}.
Ergodic Theory Dynam. Systems \textbf{22} 1 (2002), 63--97.

\bibitem{BKS} H. Bruin, A. Kaffl, D. Schleicher,  {\it Symbolic
dynamics of quadratic polynomials}. Monograph, in preparation.

\bibitem{BS} H. Bruin, D. Schleicher, {\it Symbolic dynamics of quadratic
polynomials}, Mittag Leffler Preprint {\bf 7} (2001/02). To appear as
\cite{BKS}.

\bibitem{BSnew} H. Bruin, D. Schleicher, \emph{Hausdorff dimension of 
biaccessible angles for quadratic polynomials}. Manuscript, in 
preparation.

\bibitem {DH} A. Douady, J. H. Hubbard, {\it \'Etude
dynamique des polyn\^omes complexes}, I, Publ. Math. Orsay, 1984.

\bibitem {F} K. Falconer, {\it Fractal geometry, mathematical 
foundations and applications} 2nd edn., Wiley, 2003.

\bibitem{Kiwi} J. Kiwi, \emph{Wandering orbit portraits}.
Trans. Amer. Math. Soc. \textbf{354} 4 (2002), 1473--1485.

\bibitem {M} J. Milnor, {\it Dynamics in one complex variable}, 3rd
edn., Princeton University Press, 2006.

\bibitem {Po} C. Pommerenke, {\it Boundary behavior of conformal
maps}, Springer, 1992.

\bibitem{RR} R. Radu, \emph{Hausdorff dimension and biaccessibility
for polynomial Julia sets}. Bachelor's thesis, Jacobs University, 2007.

\bibitem{SZ} D. Schleicher, S. Zakeri, {\it On biaccessible points in
the Julia set of a Cremer quadratic polynomial}. Proc. AMS {\bf 128}
3 (1999), 933--937.

\bibitem {Sm} S. Smirnov, {\it On supports of dynamical laminations
and biaccessible points in polynomial Julia sets},
Colloq. Math. \textbf{87} (2001) 2, 287--295.

\bibitem{Th} W. Thurston, \emph{On the geometry and dynamics of 
iterated rational maps}. In: D. Schleicher (ed.), \emph{Complex 
dynamics: families and friends}, AK Peters, Wellesley, MA, 2009, 
1--137.

\bibitem {Za} S. Zakeri, {\it Biaccessibility in quadratic Julia sets},
Ergod. Th. \& Dynam. Sys., \textbf{20} (2000), 1859--1883.

\bibitem {Za1} S. Zakeri, {\it External rays and the real slice of
the Mandelbrot set},
Ergod. Th. \& Dynam. Sys., \textbf{23} (2003), 637--660.


\bibitem {Zd} A. Zdunik, {\it On  biaccessible points in Julia sets
of polynomials},
Fund. Math., \textbf{163} (2000), 277--286.

\bibitem {Zd1} A. Zdunik, {\it Parabolic orbifolds and the dimension
of the maximal measure for rational maps},
Invent. Math., \textbf{99} (1990), 627--649.

\end{thebibliography}

\end{document}